\documentclass[12pt]{article}

\usepackage{comment}
\usepackage{mathtools}
\DeclarePairedDelimiter{\ceil}{\lceil}{\rceil}
\usepackage{graphicx}
\usepackage{xcolor}
\usepackage[outdir=./]{epstopdf}
\usepackage{comment}
\usepackage{caption}
\usepackage{subcaption}
\usepackage{float}
\usepackage{eurosym}
\usepackage{amsfonts}
\usepackage{fullpage}
\usepackage[margin=1in]{geometry}
\usepackage[all]{nowidow}
\usepackage{authblk}
\usepackage{textcomp}

\usepackage[numbers]{natbib}

\newcommand\Mycite[1]{%
  \citeauthor{#1}~[\citeyear{#1}]}

\begin{document}
\title{A Discrete Event Simulation Model for Coordinating Inventory Management and Material Handling in Hospitals}





\author[1a]{Amogh Bhosekar}
\affil[1]{\small Department of Industrial Engineering, Clemson University,
Freeman Hall, Clemson, SC 29634 \authorcr
  Email: abhosek@clemson.edu \authorcr
    Email: tisik@clemson.edu}

\author[2]{Sandra Ek\c{s}io\u{g}lu}
\affil[2]{Department of Industrial Engineering, University of Arkansas, 
Bell Engineering Center, Fayetteville, AR 72701 \authorcr
  Email: sandra@uark.edu} 
\author[1b]{Tu\u{g}\c{c}e I\c{s}{\i}k}
\author[3]{Robert Allen}
\affil[3]{Perioperative Services, Prisma Health,
701 Grove Rd, Greeville, SC 29605 \authorcr Email: robert.allen@prismahealth.org} 

\date{}

\maketitle

\vspace{11in}

\title{\large \bf A Discrete Event Simulation Model for Coordinating Inventory Management and Material Handling in Hospitals}

\begin{abstract}
{
For operating rooms (ORs) and hospitals, inventory management of surgical instruments and material handling decisions of perioperative services are critical to hospitals’ service levels and costs. However, efficiently integrating these decisions is challenging due to hospitals’ interdependence and the uncertainties they face. These challenges motivated the development of this study to answer the following research questions: (R1) \emph{How does the inventory level of surgical instruments, including owned, borrowed and consigned, impact the efficiency of ORs?} (R2): \emph{How do material handling activities impact the efficiency of ORs?} (R3): \emph{How do integrating decisions about inventory and material handling impact the efficiency of ORs?} Three discrete event simulation models are developed here to address these questions. Model 1, \emph{Current}, assumes no coordination of material handling and inventory decisions. Model 2, \emph{Two Batch}, assumes partial coordination, and Model 3, \emph{Just-In-Time} (JIT), assumes full coordination. These models are verified and validated using real life-data from a partnering hospital. A thorough numerical analysis indicates that, in general, coordination of inventory management of surgical instruments and material handling decisions has the potential to improve the efficiency and reduce OR costs. More specifically, a \emph{JIT} delivery of instruments used in short-duration surgeries leads to lower inventory levels without jeopardizing the service level provided.

}
\end{abstract}




%


\section{Introduction}\label{intro} 

\textbf{Motivation:} The cost of supply chain activities that support ORs contributes to a hospital’s total expenses, accounting for as much as 40\% of the operating budget in hospitals [\citenum{dobson2015configuring}]. Holding inventory of supplies and surgical instruments makes up about 10\% to 18\% of these expenses [\citenum{volland2017material}]. For example, a 2014 study conducted in acute care hospitals in California suggests that an OR costs, on average, between \$36 and \$37 each minute. The use of surgical instruments costs between \$2.50 and \$3.50 each minute [\citenum{childers2018understanding}]. Hospitals maintain large inventories to ensure that the required instruments are available for a scheduled surgery, since the lack of an instrument leads to delays. A study conducted by \Mycite{wubben2010equipment} suggests that 45.9\% of the delays in an OR happened because an instrument was unavailable [\citenum{wubben2010equipment}]. These delays resulted in longer working hours for doctors and staff, and thus, additional costs for the hospital. A surgery delayed due to a lack of instruments also negatively impacts the quality of care and adverse effects can occur [\citenum{wubben2010equipment}].

Increasing inventory levels may not necessarily eliminate these delays since some delays occur because of inefficiencies in the material handling process. For example, congestion, due to the movement of Automated Guided Vehicles (AGVs) along the narrow corridors of a hospital, can lead to delays. Based on our review of the literature, very little research evaluates the effects of inventory and material handling decisions relative to the efficiency of and service levels provided by ORs in hospitals. This gap in the literature is the main motivation for this research.


{\bf Background:} Surgeries performed in most hospitals are categorized as elective or emergency [\citenum{gupta2008appointment}]. This research focuses on elective surgeries. An elective surgery is scheduled within 12 weeks, and the exact timing of the surgery and the room assignment are finalized between 24 and 48 hours before the day of the surgery. The OR scheduler makes these assignments after considering the availability and preferences of the surgeon and surgical staff, as well as the availability of the required equipment and instruments. These assignments impact the availability of instruments for the rest of the day. Emergency surgeries are incorporated into the daily schedule since the timing of intervention is critical for patients’ safety [\citenum{gupta2008appointment}]. To maintain high service levels, some hospitals develop plans to ensure that sufficient ORs, instruments, and equipment are available. The limited availability of instruments restricts hospitals from utilizing ORs and surgeons’ time efficiently. In most hospitals, an instrument is not used more than once the same day because it should be decontaminated and sterilized before reuse. This process can take up to 3 or 4 hours. Most of the time, an OR operates 12 hours each day, and most of the surgeries last no more than 5 hours.

\textbf{Research Questions:} Surgical instruments are categorized as: (1) owned by the hospital, (2) borrowed from other hospitals or doctors, or (3) consigned by a vendor who owns the instrument [\citenum{chobin2015best}]. The process of adopting a new surgical instrument is initiated upon a surgeon’s request. Next, the hospital evaluates whether buying, renting, or consigning this instrument is the best option. The main factors impacting this decision are the selling price of the instrument and the frequency of its use. Typically, a hospital would not purchase an instrument if only used in rare or specialty surgeries [\citenum{chobin2015best}]. In such a case, the hospital would consign or borrow the instrument and pay the owner upon its use. A hospital has several other reasons to borrow instruments, such as to match demand and supply, accommodate doctors’ requests for scheduling consecutive surgeries during a given day, continue operations on a limited budget, or mitigate a lack of storage space [\citenum{seavey2010reducing}]. According to \Mycite{seavey2010reducing}, such practices lead to inefficiencies because borrowed instruments create an additional workload for the sterile processing department (SPD) since the hospital is required to maintain documentation and pack and sterilize instruments [\citenum{seavey2010reducing}]. Additionally, some instruments have special cleaning procedures, which may differ from other procedures for instruments owned by the hospital. Following these procedures adds to employees’ workloads. Also, consigned instruments stored at the hospital occupy additional storage space. A recent study conducted in a major US academic hospital suggests that half of the instruments are consigned, and their cost is, on average, 12\% more than instruments owned by the hospital [\citenum{mandava2017how}]. These challenges motivate the first research question: \textbf{(R1)} \emph{How does the inventory level of surgical instruments, including owned, borrowed, and consigned, impact the efficiency of ORs?} Typically, the inventory level is determined by the total number of surgeries scheduled in a day, the daily schedule of surgeries that use the same instrument, the processing capacity of the central sterile storage division (CSSD), and the schedule of material handling activities. This study presents data-driven, discrete event simulation (DES) models and a numerical study that evaluates how inventory levels impact the utilization of instruments and the delay of surgeries. The models are validated with data from a partnering US-based hospital.

Instruments are delivered to ORs via containers called case carts that are transported by carriers, such as staff or AGVs. After a surgery, soiled instruments are delivered to the CSSD using the same carriers. Inefficiencies in material handling activities lead to delays, which impact the availability of instruments. Additionally, the duration of a surgery is uncertain; thus, a surgery may take longer than planned, keeping instruments unavailable. For these reasons, some hospitals deliver instruments to a storage area beside the OR the night before the surgery. Such a practice ensures that the instruments required are available during the surgery. As a result, the same instrument cannot be reused in other surgeries scheduled on the same day. Alternatively, an instrument could be delivered to an OR directly from the CSSD a short time before the start of a surgery, the JIT delivery approach. Such a practice increases the utilization of instruments used in short-duration surgeries performed earlier in the day. This approach can also lead to lower inventory levels and lower inventory holding costs. However, such an approach requires coordination of material handling, instrument decontamination and sterilization, and OR scheduling. These challenges motivate the second research question: \textbf{(R2)} \emph{How do material handling activities impact the efficiency of ORs?} A numerical analysis via the proposed DES models is conducted to answer this question. These models take different approaches to the handling of instruments. Each material handling approach follows a different schedule of delivering case carts to ORs. For each approach, a numerical study is conducted to evaluate how the number of carriers impacts travel time, congestion, the utilization of carriers, delivery time, and the utilization of instruments.

Coordinating decisions about material handling schedules and instrument inventory is challenging. For example, the decision to reduce the inventory of an instrument limits the time that instrument is available. This, in turn, negatively impacts the flexibility of scheduling a surgery needing the instrument and, therefore, the service level provided. The problem becomes even more challenging when other important considerations are added, such as the inconsistent material handling system, stochastic demand, and uncertain surgery duration. These challenges motivate the following research question: \textbf{(R3)} \emph{How do integrating decisions about inventory and material handling impact the efficiency of ORs?} To answer this question, another numerical study using the DES models is conducted.

\textbf{Contributions:} The proposed research offers several important contributions: (\textit{i}) This study highlights the role of coordinating decisions between material handling and inventory management in improving OR efficiency and reducing costs while maintaining high service levels. In particular, this work demonstrates that JIT delivery of surgical cases for short-duration surgeries can potentially improve the efficiency of ORs and reduce the cost of healthcare. (\textit{ii}) This study develops a real-life case study using data from a US-based hospital. The proposed material handling approaches, which are intuitive and easy to implement, are verified and validated using historical data. The results of the proposed analysis have inspired the partner hospital featured in this study to make improvements in material handling and inventory management practices. While other healthcare facilities may not choose to implement the models presented here, they can learn from these practices.

\textbf{Outline:} The rest of this article is organized as follows: Section 2 reviews the literature relevant to this work, and Section 3 provides a detailed description of the problem. Section 4 describes the proposed simulation models, and Section 5 introduces the case study. Section 6 discusses the results of the experiments, and finally, Section 7 summarizes the results and presents concluding remarks.

\section{Literature Review}
The main stream of existing literature relevant to this work is inventory management of reusable surgical instruments. Since AGVs are used as carriers by this study's partner hospital, as well as in many others, the literature that discusses the use of AGV systems for material handling in hospitals is also reviewed.

{\bf Inventory Management of Reusable Surgical Instruments:}

The cost of ORs is impacted by the availability of surgical supplies and implants [\citenum{chasseigne2018assessing}]. Surgical supplies include soft goods and instruments required for surgery. These supplies can be either reusable or disposable. Numerous studies show that the cost of reusable supplies is significantly lower than the cost of disposable supplies [\citenum{demoulin1996cost}, \citenum{eddie1996comparison}, \citenum{schaer1995single}, \citenum{manatakis2014reducing}, \citenum{adler2005comparison}], {and} disposable supplies negatively impact the environment [\citenum{adler2005comparison}]. \Mycite{chasseigne2018assessing} conduct a study to evaluate the cost of opened, unused soft goods and instruments in a French hospital. They reported that wasted supplies have a median cost of \euro 4.1 per procedure, which accounts for about 20.1\% of the cost of surgical supplies. However, most hospitals do not have standardized procedures to manage the inventory of surgical supplies [\citenum{ahmadi2018inventory}]. In their review paper, \Mycite{ahmadi2018inventory} indicate that inventory management of sterile instruments requires three important considerations: instrument and quantity assignment for each tray-type, the tray-type’s assignment to a surgeon or procedure, and the number of trays carried by the hospital. Decisions related to the first two considerations are impacted by the surgeon’s preferences, indicated in the doctor preference card (DPC).

The cost of surgical supplies can be reduced in several ways, including by 1) improving the accuracy of the DPCs, 2) increasing surgeon awareness, and 3) standardizing surgical techniques. Accuracy of the DPC can be improved by reviewing it periodically [\citenum{harvey2017physician}, \citenum{farrelly2017surgical}] or by recording which instruments are used on a tray and removing the instruments that are not used [\citenum{nast2019decreasing}, \citenum{dyas2018reducing}]. For example, \Mycite{harvey2017physician} show that engaging physicians in the review of the corresponding DPC led to the removal of 109 disposable supplies and the elimination of 3 reusable instrument trays. Consequently, the cost of a case cart was reduced by \$16 on average. According to a survey conducted by \Mycite{jackson2016surgeon}, surgeons often underestimate the cost of expensive items and overestimate the cost of less expensive items due to internal bias and cost ignorance [\citenum{jackson2016surgeon}]. Thus, the cost of a surgical procedure can be reduced by increasing surgeons’ awareness of standardized operating equipment and the cost of instruments [\citenum{gitelis2015educating}, \citenum{avansino2013standardization}]. Finally, work by \Mycite{skarda2015one} shows that standardization of surgical techniques can significantly reduce operating costs without impacting the quality of a procedure [\citenum{skarda2015one}].

\Mycite{stockert2014assessing} indicate that tailored and streamlined tray compositions lead to significant cost savings [\citenum{stockert2014assessing}]. Additionally, surgeons prefer trays that have fewer unsolicited instruments [\citenum{dobson2015configuring}, \citenum{stockert2014assessing}]. Several optimization models have been developed to solve the tray optimization problem and address tray composition and inventory management for reusable surgical instruments. The objective of this problem is to minimize an OR's cost by optimizing the number of trays utilized and the amount of inventory supplied. The problem also addresses surgeon preferences for instruments. \Mycite{dobson2015configuring} develop a linear integer programming formulation and propose a heuristic algorithm to obtain a solution to this problem [\citenum{dobson2015configuring}]. \Mycite{reymondon2008optimization} propose a resource sharing method for reusable devices [\citenum{reymondon2008optimization}]. The objective is to minimize storage, processing, and wastage costs for supplies that have not been used. \Mycite{van2008optimizing} propose a deterministic model that minimizes the storage and delivery cost of instruments by optimizing tray composition [\citenum{van2008optimizing}]. \Mycite{ahmadi2019bi} present a bi-objective optimization model for configuration of surgical trays with ergonomic considerations [\citenum{ahmadi2019bi}]. The first objective function minimizes the total number of assembled tray types, and the second objective function minimizes the total number of instruments that were not requested. They use the $\epsilon$- constraint method to obtain the Pareto-optimal front. \Mycite{dollevoet2018solution} develop an exact integer linear programming formulation, a row and column generation approach, a greedy heuristic, and some metaheuristics. These approaches are evaluated based on the average computation time, the average value of the objective function, and the number of solutions for which optimality is proven.

{\bf AGV Systems and Operations in Hospitals:}

The existing literature on AGV systems focuses on fleet size selection. Simulation and optimization models are proposed to identify AGV fleet size (\Mycite{choobineh2012fleet}, \Mycite{maxwell1982design}, \Mycite{arifin2000determination}, \Mycite{rajotia1998determination}, \Mycite{sinriech1992economic}, \Mycite{egbelu1987use}, \Mycite{tanchoco1987determination}, and \Mycite{Bhosekar:2018}). 


\Mycite{katevas2001mobile} discusses several factors that must be considered to design a mobile robotic system for healthcare applications, and his work provides several guidelines for researchers to improve these designs. \Mycite{rossetti2000simulation} compare the manual delivery of clinical and pharmaceutical items with the performance of robotic courier delivery [\citenum{rossetti2000simulation}]. They use cost, turnaround time, variability of turnaround time, cycle time, and utilization as performance measures. The proposed simulation model shows that using robotic delivery is economically viable and improves the performance measures listed above. \Mycite{rossetti2001multi} use Analytic Hierarchy Process to build a decision problem that evaluates the performance of a robotic healthcare delivery system based on technical, economical, and several other factors. Their proposed simulation model assesses the technical {factors} that include speed of robot and human couriers based on the arrival rates of visitors who request the elevator, elevator availability, the arrival rates of the delivery items that request robots, and robot availability. In their case study, \Mycite{chikul2017technology} compare three supply chain models that use: a) manual inventory check and delivery, b) RFID inventory check and manual delivery, and c) manual inventory check and AGV-based material handling. This study shows that combined RFID tracking and AGV-based delivery maximizes cost savings in the supply chain model and yields ergonomic benefits due to reduced manpower requirements. Via a simulation-based case study, \Mycite{pedan2017implementation} identify potential benefits of utilizing an AGV system in a hospital. Finally, \Mycite{fragapane2018material} evaluate the impact that material and information flow have on costs at a Norway hospital.

{Different from this literature which focuses either on improving the management of inventory of surgical instruments, or improving material handling activities in hospitals, our proposed work focuses on evaluating the impacts of integrating material handling and inventory management decisions. Via our  numerical analysis we show that coordinating these decisions leads to reduced inventory levels, reduced  number of AGVs used and increased utilization of AGVs, while the quality of service provided remains intact. }

\section{Problem Description}\label{mhp}
\begin{figure}[ht]\centering{
{\includegraphics[width=0.75\textwidth]{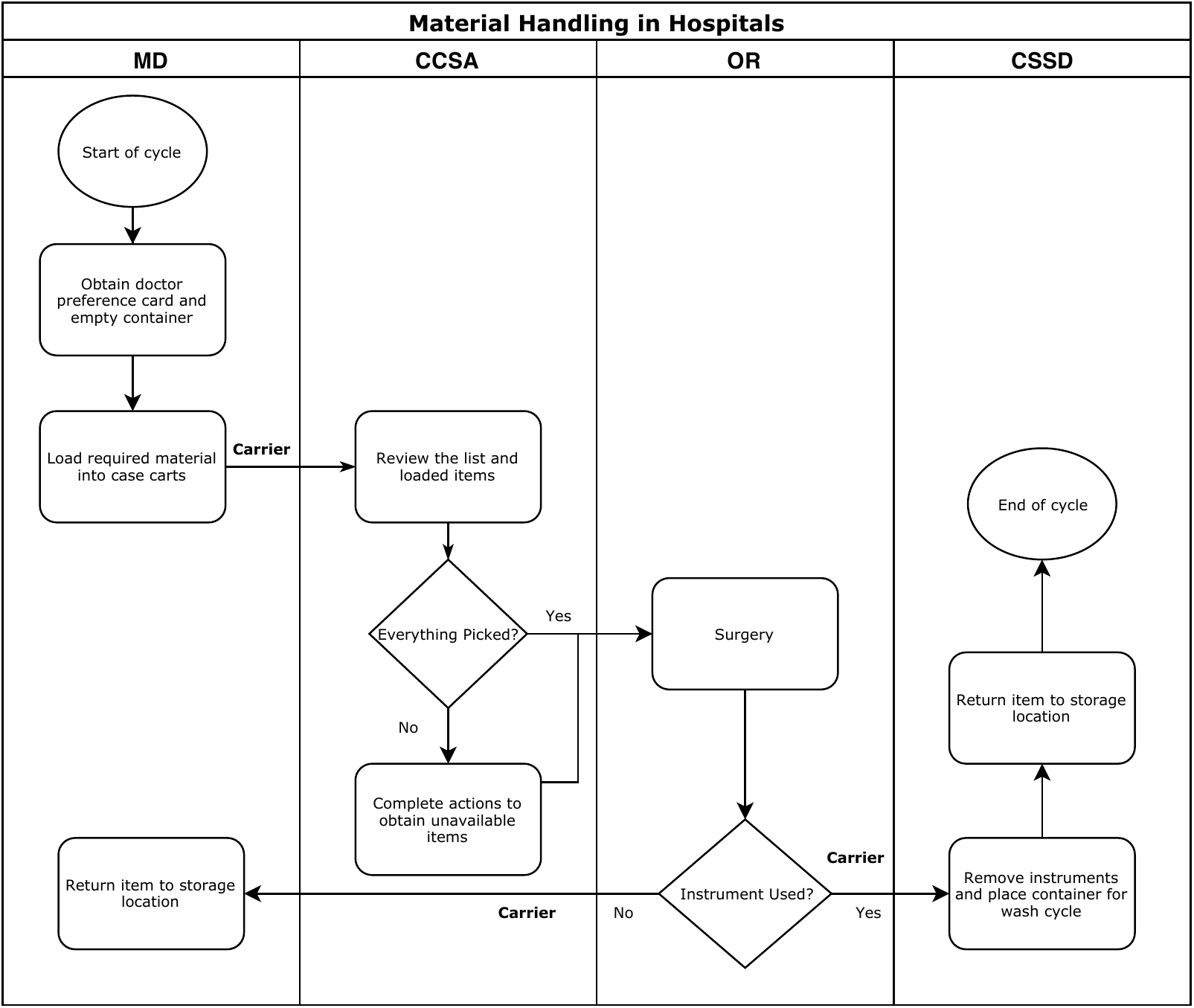}}
\caption{Material Handling Process \label{fig: MHProcess}}
}
\end{figure}

\begin{figure}[ht]
{\includegraphics[width=\textwidth]{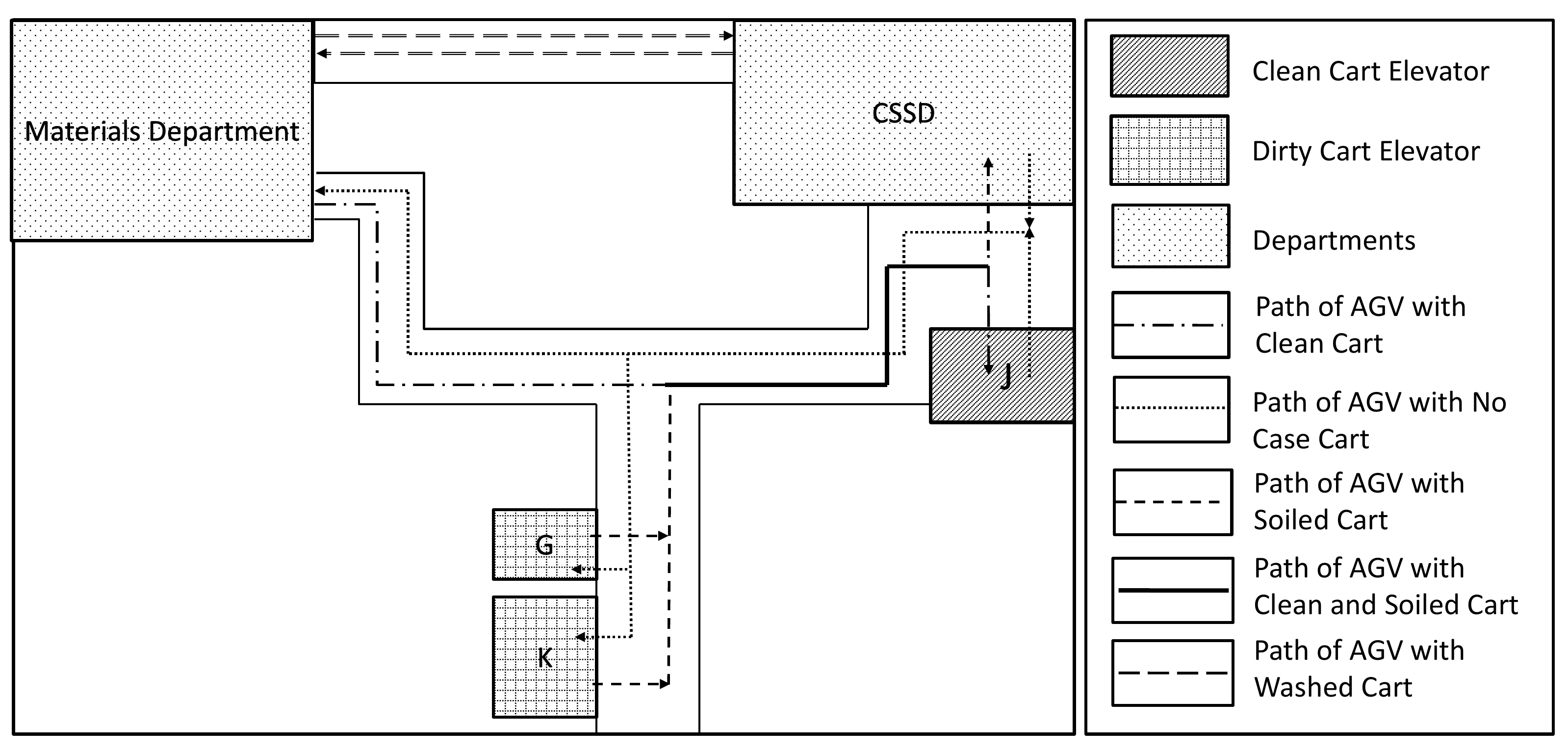}}
\caption{Map - GMH Floor Map}
{\label{fig: Model1}}
\end{figure}

\textbf{Material Handling and Inventory Management Processes:} Figure \ref{fig: MHProcess} describes a typical material handling process for the delivery of surgical case carts in a hospital. The process begins by creating a detailed schedule of surgeries. This schedule is prepared by the OR manager. Based on the schedule and the doctors’ preferences, a list of instruments and soft goods is prepared and submitted to the materials division (MD). For each surgery, a clean case cart is loaded with the requested instruments, soft goods, and implants. These case carts are moved to pick-up/drop-off stations for carriers to pick up. Then, the clean case carts are moved from the MD to the case cart storage area (CCSA). At the CCSA, each case cart is inspected to ensure that it contains the required materials. The case carts are held at the CCSA until they are moved to the corresponding OR at the time of the surgery. The case carts are delivered to ORs prior to the surgeries. ORs are divided into separate cores based on the specialties they serve. Specialty instruments and implants required for the surgical cases, which are stored in the OR cores, are added to the case carts before the surgery. After the surgery, the instruments and case carts are considered soiled and should be decontaminated. The soiled carts and instruments are transported to the CSSD by the carrier. The instruments and case carts are washed and sterilized at the CSSD. The specialty surgical instruments are returned to the corresponding OR cores. This process ensures the availability of instruments before the scheduled surgery. Figure \ref{fig: Model1} presents the locations of the departments and paths traversed by the carriers (i.e., AGVs) at our partnering hospital. 

\begin{table}[ht]
\centering
\caption{AGV Movements by the Time of the Day \label{Table: ByTime}}
{\begin{tabular}{llllll}
\hline
Route & Time & No. of & \multicolumn{2}{c}{Travel Time [Min]} & Coefficient of \\\cline{4-5}
 & Interval & Trips & Average & Std. Dev.  & Variation \\ \hline
 & 12 am-3 am & 131 & 4.66 & 10.26 & 2.2\\
 & 3 am-6 am & 227 & 5.96 & 9.55 & 1.6 \\
 & 6 am-9 am & 112 & 5.27 & 6.17 & 1.17 \\
Materials Department - & 9 am-12 pm & 80 & 5.8 & 2.4 & 0.41 \\
Case Cart Storage Area & 12 pm-3 pm & 101 & 5.33 & 3.69 & 0.69 \\
 & 3 pm-7 pm & \textbf{1,416} & \textbf{8.94} & \textbf{6.49} & 0.73 \\
 & 7 pm-9 pm & 254 & 5.67 & 4.45 & 0.78\\
 & 9 pm-12 am & 196 & 4.88 & 5.72 & 1.17\\ \hline
  & 12 am-3 am & 44 & 7.28 & 6.16 & 0.85\\
 & 3 am-6 am & 53 & 6.87 & 6.16 & 0.9\\
 & 6 am-9 am & 146 & 5.55 & 3.27 & 0.59\\
2nd Floor Soiled & 9 am-12 pm & 981 & 4.56 & 4.04 & 0.88 \\
Cart Storage - CSSD & 12 pm-3 pm & 882 & 5.17 & 2.56 & 0.49\\
 & 3 pm-7 pm & 753 & \textbf{9.71} & \textbf{8.51} & 0.88 \\
 & 7 pm-9 pm & 126 & 6.65 & 3.04 & 0.46\\
 & 9 pm-12 am & 77 & 5.45 & 1.01 & 0.19\\ \hline
\end{tabular}}\\
{\footnotesize{This table was obtained from the prior research. [\citenum{Bhosekar:2018}]}}
\end{table}

To reduce the cost of inventory, hospitals need to coordinate inventory management and material handling decisions. This coordination becomes ever more important in face of uncertainty. For example, if surgery duration and travel time of carriers is fixed, hospitals can calculate the necessary inventory levels with certainty and decide how many instruments to loan or consign. However, in order to ensure high service level under uncertainty, many hospitals keep large inventories, loan instruments, and prepare/deliver case carts one day before the surgery. We observed at our partner hospital that, if the delivery of instruments from the CSSD to ORs was completed within a short time before the start of the surgery, the instruments could be reused within the same day. Such an approach has the potential to lead to a reduction in the cost of using loaned or consigned instruments. This observation led the development of the material handling process proposed next.

\begin{table}[ht]
\centering
\caption{A List of Instruments Used \label{Table: IAGMH}} 
{\begin{tabular}{r|rr}
\hline
\multicolumn{1}{c|}{\bf Type} & {\bf Total Number} & {\bf Total in \%} 
\\ \hline
Loaner & 266 & 5\% \\
Consigned & 1,095 & 19\% \\
Owned & 3,507 & 61\% \\
Other Services & 927 & 16\% \\ \hline
{\bf Total} & 5,795 & 100\% \\ \hline
\end{tabular}}
{}
\end{table}

\textbf{Experimental Setup:} The movement of clean instruments to ORs and the movement of soiled instruments to the CSSD affect inventory {availability} and the starting times of surgeries. This is the reason why some hospitals, like our partner hospital, prepare and deliver the surgical case carts to the CCSA one day in advance. Such a practice ensures the availability of surgical instruments, but a number of inefficiencies results regarding material handling and inventory management. For example, Table \ref{Table: ByTime} summarizes the data obtained on the travel times of AGVs during different times of the day at the partner hospital. The data shows that the average travel time  and the corresponding standard deviation are highest during 3 pm to 7 pm. This is because the clean surgical carts are being delivered to the CSSA for elective surgeries. The consequent increase in the number of AGV movements leads to congestion and thus, longer travel times for every AGV that uses the same path. These delays lead to an increased inventory of instruments since an instrument cannot be reused in different surgeries scheduled on the same day. Table \ref{Table: IAGMH} lists the types of the total number of instruments used by the partner hospital and the percentage of each type. Notice that about 24\% of the instruments used are either loaned or consigned.

This study evaluates three approaches to delivering surgical supplies to ORs and compares their performances. The first approach, Model 1, is referred to as the \emph{Current} approach and assumes that materials required by surgeries are delivered to the CCSA the night before the surgery. The \emph{Current} approach is the ongoing practice of the partner hospital in the data presented here.

Next, the \emph{Two Batch} approach, Model 2, assumes that materials required by surgeries scheduled in the morning are delivered to the CCSA the previous evening, and the materials required by surgeries scheduled in the afternoon are delivered in the morning on the day of the surgery. This approach provides the opportunity to reuse the instruments from the surgeries scheduled later in the day. Since the CSSD works 24 hours each day, the instruments can be washed overnight and delivered in the morning. Since instruments are delivered a few hours in advance of the surgery, the staff has an abundant amount of time to intervene if an instrument becomes unavailable. Thus, the risk of instruments not being delivered on time is only minimal and does not affect the quality of care in the hospital.

Finally, the \emph{Just-in-Time} approach, Model 3, assumes that materials required are delivered shortly before the start of the surgery. The time between the delivery of surgical supplies and the surgery, referred to as the delivery interval, needs to be determined and impacts the inventory levels. The inventory level required increases with the delivery interval. For example, consider two surgeries that require the same instrument and are scheduled on the same day. In the current system, a hospital must have two sets of the identical instruments since they are delivered to the CCSA the day before the surgery. If the delivery interval is 1 hour, the instrument can be sterilized and delivered before the subsequent surgery, provided that the two surgeries are scheduled several hours apart. In this case, the hospital needs only one instrument. However, if the delivery interval is chosen to be less than 3 hours, then, to avoid any delays of the second surgery, two instruments are needed. Note that the implementation of JIT and other lean methods in healthcare, unlike with manufacturing, should be considered with caution because such practices could delay surgeries and jeopardize the well being of patients.

{Three performance measures are used to compare the proposed approaches: (\textit{i}) the average delay of a surgery's start time, which is a measure of the service level provided, (\textit{ii}) the number of instruments inventoried, which  measures the efficiency of the inventory system, and (\textit{iii}) the number of carriers required for each proposed material handling approach, which measure the efficiency of the material handling systems.}



The delays of a surgery start time are separated into two categories: delays due to carriers, e.g., long travel time because of congestion or unavailability of carriers, and delays due to unavailable instruments. Delays due to carriers can create challenges for the JIT approach, but these delays can be reduced by optimizing fleet size. Delays due to unavailable instruments are caused either by delays in the delivery of soiled case carts or by an increase in the number of emergency surgeries. A delay in the delivery of soiled case carts subsequently delays the cleaning process of instruments and carts, which leads to the delay of the start of the next surgery that uses the same instrument. Delays due to unavailable instruments can be reduced by optimizing the inventory level. In this research, simulation experiments are conducted to determine the optimal fleet size and inventory level under each proposed material handling approach. Based on the results of these experiments, the delivery interval that optimizes the performance measures identified is also determined.

Each of the proposed material handling approaches requires a different number of carriers to deliver materials on time. This number is impacted by the surgery schedule and the material handling process. For instance, the number of carriers needed for the JIT delivery approach is lowest since the delivery of case carts is spread throughout the day. The number of carriers needed by the \emph{Current} delivery approach is larger because the delivery of case carts is completed within a short time period. The number of carriers needed also depends on the total number of cases scheduled and the spread of the schedule. A tight schedule would require more carriers to complete material handling on time.

{\bf Limitations of this Research:}

{\bf Model:} The research proposed here is conducted in collaboration with Greenville Memorial Hospital (GMH), a US-based hospital located in Greenville, South Carolina. The models presented here are  motivated by the material handling and inventory management practices at GMH. The research team worked closely with the perioperative services department, which consists of the MD, the CSSD, and the OR Division. GMH uses AGVs to transport surgical case carts to and from ORs. The problem setting proposed here and the assumptions made are influenced by the practices at GMH. The models presented here are a valuable contribution to the literature because, based on a careful review of the literature, similar practices are followed by other hospitals for material handling and inventory management of surgical instruments. 

\noindent {\bf Data:} Nine months of real-life data are used to develop the case study. This data includes information about the number of surgical cases each day and ranges over a time period long enough to observe how seasonality impacts the number of surgical cases. Ideally, larger amounts of data would be available, but that does not apply here. 

\section{Simulation Model}

DES models are developed to evaluate and compare the three approaches proposed for the delivery of surgical supplies. These models are created in ARENA simulation software by Rockwell Automation. An entity type represents a surgery type, and each entity represents a surgical case of a particular type. An entity has three attributes: \emph{duration, starting time, and type}. \emph{Duration} is randomly generated using the distributions listed in Table 3. These distributions are derived from the data collected at GMH. The \emph{starting time} and \emph{type} are fed to the model from the actual data. Other entities are used to control the movement of AGVs and elevators, as well as to handle other specific requirements, such as calculating the value of certain variables (e.g., the number of AGVs to activate each day). ORs, case carts, cart washers, and elevators were modeled as resources. Variables are used to track the number of busy resources. A guided path transporter network is developed with intersections and links to replicate the movement of AGVs along the corridors of the hospital. This network was constructed using actual distances obtained from a GMH floor map. The links of the network are unidirectional, bidirectional, or spurs (dead ends). The intersections represent the areas where two or more links intersect. The intersections allow AGVs to make turns and move from one link to the next, following their routes. Intersections are also used to represent pick-up/drop-off stations. A spur link marks the end of a route. Departments can only handle a certain number of AGVs, and their processing capacity is limited by variables.

\begin{table}[ht]
\caption{Input Parameters: Surgery Duration \label{Table: LOS}}
{\footnotesize \begin{tabular}{l|cc|l|l|c}
\hline
Service & From & To & Distribution & Expression (Length of Surgery) & Squared Error \\ \hline
ENT Surgery& 00:00 & 08:00 & Lognormal & LOGN(2.02, 2.12) & 0.008 \\
 & 08:00 & 14:00 & Lognormal & LOGN(1.62, 1.2) & 0.007 \\
 & 14:00 & 00:00 & Lognormal & LOGN(1.23, 0.672) & 0.003 \\ \hline
Gynecology Service & 07:00 & 08:00 & Beta & 0.01 + 4.81 * BETA(2.85, 4.03) & 0.009 \\
 & 15:00 & 16:00 & Lognormal & 0.27 + LOGN(0.965, 0.511) & 0.005 \\
 & 16:00 & 07:00 & Lognormal & LOGN(1.65, 0.859) & 0.011 \\ \hline
Neurological Surgery  & 00:00 & 09:00 & Gamma & GAMM(0.494, 5.44) & 0.007 \\
 & 09:00 & 13:00 & Erlang & ERLA(0.454, 5) & 0.005 \\
 & 13:00 & 00:00 & Beta & 12 * BETA(4.95, 25.8) & 0.028 \\ \hline
Ortho Trauma Surgery  & 0:00 & 8:00 & Erlang & ERLA(0.587, 5) & 0.002 \\
 & 08:00 & 14:00 & Lognormal & LOGN(2.57, 1.29) & 0.004 \\
 & 14:00 & 00:00 & Lognormal & LOGN(2.09, 0.983) & 0.004 \\ \hline
Pediatric Surgery & 00:00 & 00:00 & Lognormal & LOGN(1.35, 0.693) & 0.011 \\ \hline
Urology Surgery & 00:00 & 07:00 & Lognormal & LOGN(1.63, 0.975) & 0.001 \\
 & 07:00 & 08:00 & Lognormal & LOGN(1.2, 0.821) & 0.007 \\
 & 08:00 & 00:00 & Erlang & ERLA(0.244, 4) & 0.011 \\ \hline
Vascular Surgery  & 00:00 & 07:00 & Beta & 0.03 + 8.97 * BETA(0.97, 1.78) & 0.004 \\
 & 07:00 & 09:00 & Gamma & GAMM(0.608, 3.58) & 0.017 \\
 & 09:00 & 14:00 & Gamma & GAMM(0.42, 4.39) & 0.025 \\
 & 14:00 & 00:00 & Triangular & TRIA(0.13, 0.83, 3.54) & 0.011 \\ \hline
\end{tabular}}
{}
\end{table}

\begin{figure}[ht]
{\includegraphics[width=\textwidth]{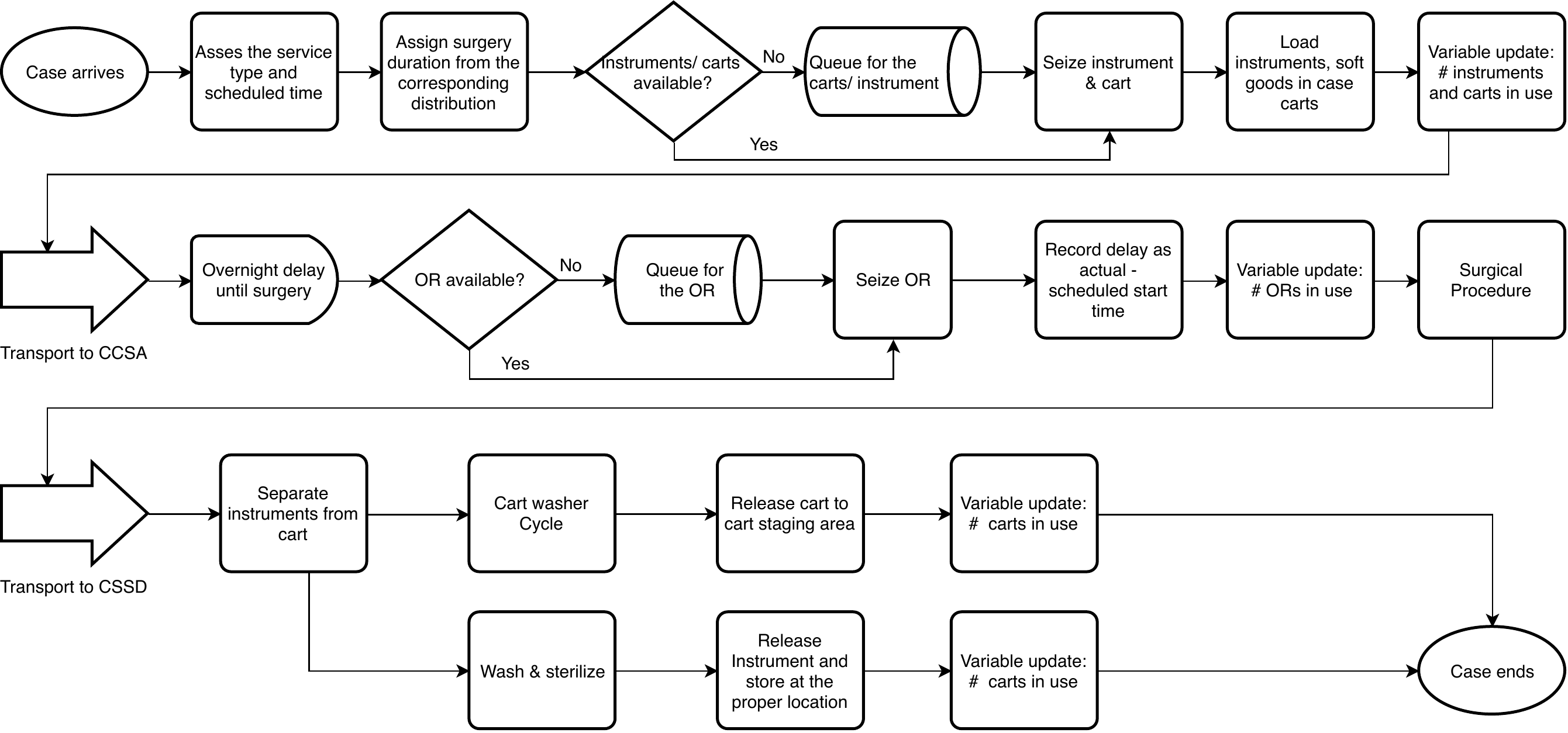}}
\caption{Flowchart of the Simulation Model \label{fig: SimulationFlowchart}}
{}
\end{figure}

The first DES model, Model 1, depicts the \emph{Current} material handling approach used at GMH. Figure \ref{fig: SimulationFlowchart} describes this model and Table \ref{tab:param} summarizes the values of its input parameters. In this model, the release of entities begins at 3 pm. The start time of these entities takes place after 6 am the next day. Next, the availability of instruments is checked using the decide module. If an instrument is not available, the case cart is held at MD until the instrument becomes available. An available instrument seizes a case cart and is delivered to the CCSA. There, the entity is held until the scheduled start time of the surgery. At this point in time, the entity seizes an available OR for the duration of the surgery. At the end of a surgery, the OR is released, and the corresponding case cart and instrument are moved to the CSSD to wash and sterilize. The resources at the CSSD are seized for the duration of service. The variables which record the number of busy units are updated when resources are released.

The second DES model, Model 2, depicts the \emph{Two-Batch} material handling approach. In this model, entities are released twice a day, at 6 am and 3 pm. The entities released at 6 am have a start time between after 12 pm the same day. Next, these entities follow a similar procedure as described above in Model 1. Entities released at 3 pm have a start time between 6 am and noon the next day. These entities are held until the next morning using the hold module, and then, they follow the procedure outlined above.

The third DES model, Model 3, depicts the \emph{JIT} approach. In this model, entities are released one hour prior to their \emph{start time}. This delivery interval was chosen based on the results obtained in Table \ref{Table:CCS}. Next, these entities follow the same procedure outlined above. In every model, the delivery of soiled case carts begins as soon as a surgery is completed.


\begin{table}[h]
\caption{Summary of Input Parameters }
\label{tab:param}
\resizebox{\textwidth}{!}{%
\begin{tabular}{lll}
\hline
\textbf{Parameter} & \textbf{Source} & \textbf{Description} \\ \hline
Entity Creation Time & Surgery schedule data & Read from the data \\
Attribute Duration & Surgery schedule data & Random variable from \\
 &  & the corresponding distribution \\
Network link distances & GMH floor maps & Read from the data \\
No. of Case carts & AGV system data & 110 \\
No. of ORs & GMH Survey & 32 \\
No. of loading personnel & GMH Survey & 4 \\
No. of AGVs & AGV system data & {[}6,8,10{]} \\
Capacity of elevators & GMH Survey & {[}2,2{]} \\
Capacity of cart washers & GMH Survey & 3 \\
Cart loading delay & GMH Survey & Triangular(2,3,5) minutes \\
Cart washing delay & GMH Survey &  20 Minutes \\
Elevator movement delay to carry AGV & GMH Survey &  40 seconds \\
Cart loading unloading delay & GMH Survey &  15 seconds \\
Instrument washing delay & GMH Survey &  3 hours \\ \hline
\end{tabular}%
}
\end{table}

\section{A Case Study}
{\bf Input Data Analysis:}
The main objective of the data collection and analysis is to evaluate the impacts that the \emph{Current} material handling approach has on the inventory of surgical instruments. The data collected is used to develop the DES models.

The data is presented in the following sets: The first data set provides information about the surgeries scheduled at GMH from Jan. 1, 2018, until Sept. 11, 2018. This data includes the surgery identification number (ID), OR’s ID , the date of the surgery, the scheduled start and finish times of a surgery, the type of surgery (i.e., vascular, orthopedic, neurological etc.), information about the surgeon, the primary procedure, and the instruments requested. The second data set provides information about the surgical instruments used at GMH. This data presents the instrument ID, the type of surgery the instrument is used for, the inventory level, and information regarding its ownership.

The hospital offers 46 different surgical services. Our experimental analysis focuses on the following 7 types of surgeries: ENT, pediatric, ortho trauma, neurology, gynecology, urology, and vascular. We focus on these surgical services because they are scheduled multiple times each day. Therefore, there is an opportunity to reduce the size and cost of inventory by reusing some of the instruments. The duration of a surgery is calculated using the actual start and finish times.

Surgeries are grouped based on service type, duration, and scheduled start times. For each service type, an hypothesis test is conducted to evaluate whether the duration of surgeries within each service type differ based on the starting time of the given surgery. When differences were observed, the distribution of surgery duration was estimated separately. Otherwise, the data was used to derive a single distribution for surgeries of the same type that were started at different times of the day. The results of the hypothesis test generated the input parameters used in the simulation model. For example, the surgery duration differs based on the time of the day the surgery is scheduled, by day of the week, and also by service type. A continuous distribution was fitted using the Input Analyzer of Rockwell Automation to represent the surgery duration. Table \ref{Table: LOS} shows the service types, distribution of the length of surgeries, and the squared error. The real-life scheduled start times of the surgeries are used in the simulation model obtained from the data set and presented here.

Table \ref{Table: NOS} summarizes the total number of surgical cases scheduled between Jan. 1, 2018, and Sept. 11, 2018. Here, only the surgery types that were scheduled more frequently are listed. Each of these surgery type is scheduled more than once a day and requires multiple instruments of the same kind. For each surgery type, only one set of instruments, common to all the surgical cases of that type, is used. Table \ref{Table: IISI} lists the instruments selected for this study and their corresponding inventory levels.

GMH carries multiple instruments for each surgery type for three main reasons: First, the same surgery could be scheduled more than once in the same day if the hospital follows a block schedule approach. This approach assigns the same block of time to a surgeon or a group of surgeons who perform similar procedures every week because surgeons perform back-to-back specialty surgeries in the assigned blocks and use similar instruments. The Current material handling system requires that every instrument is available one day before the surgery. Second, surgeons of different specialties may request the same instrument for the same procedure. Third, the hospital carries safety stock to respond to instrument-related incidents, such as dropping or breaking an instrument during a surgical procedures.

\begin{table}[ht]
\caption{Input Parameters: Number of Surgeries \label{Table: NOS}}
\resizebox{\textwidth}{!}{
{\begin{tabular}{l|ccccccc|c}
\hline
\textbf{Service} & \textbf{Sunday} & \textbf{Monday} & \textbf{Tuesday} & \textbf{Wednesday} & \textbf{Thursday} & \textbf{Friday} & \textbf{Saturday} & \textbf{Total} \\ \hline
\begin{tabular}[c]{@{}c@{}} ENT  Surgery\end{tabular} & 25 & 295 & 148 & 264 & 231 & 302 & 20 & 1,285 \\
\begin{tabular}[c]{@{}c@{}}Gynecology  Service\end{tabular} & 17 & 227 & 133 & 181 & 176 & 198 & 13 & 945 \\
\begin{tabular}[c]{@{}c@{}}Neurological  Surgery\end{tabular} & 22 & 174 & 168 & 293 & 163 & 225 & 17 & 1,062 \\
\begin{tabular}[c]{@{}c@{}}Ortho Trauma Surgery\end{tabular} & 2 & 205 & 171 & 174 & 206 & 207 & 56 & 1,021 \\
\begin{tabular}[c]{@{}c@{}}Pediatric  Surgery\end{tabular} & 62 & 145 & 248 & 153 & 276 & 158 & 79 & 1,121 \\
\begin{tabular}[c]{@{}c@{}}Urology   Surgery\end{tabular} & 39 & 293 & 333 & 298 & 382 & 466 & 72 & 1,883 \\
\begin{tabular}[c]{@{}c@{}}Vascular   Surgery\end{tabular} & 61 & 141 & 242 & 224 & 241 & 235 & 81 & 1,225 \\ \hline
{\bf Total} & 228 & 1,480 & 1,443 & 1,587 & 1,675 & 1,791 & 338 & 8,542 \\ \hline
\end{tabular}}
{}
}
\end{table}

\begin{table}[ht]
\centering
\caption{Number of Instruments in the Inventory \label{Table: IISI}}
{\begin{tabular}{llc}
\hline
\textbf{Service Type} & \textbf{Instrument} & \textbf{Inventory} \\ \hline
{ENT Surgery} & Set T \& A  GMMC  1047 & 10 \\
{Gynecology Service} & Set D \& C mini  GMMC  15896 & 10 \\
{Neurological Surgery} & Set Back Neuro  GMMC  1341 & 12 \\
{Ortho Trauma Surgery} & Set Minor Ortho  GMMC  100031 & 17 \\
{Pediatric Surgery} & Set Pediatric Minor  GMMC  1247 & 8 \\
{Urology Surgery} & Ureteroscope 7.5 Comp  GMMC  12656 & 18 \\
{Vascular Surgery} & Probe Doppler Pencil 8.1 GMMC  1824 & 25 \\ \hline
\end{tabular}}
{}
\end{table}

{\bf Verification and Validation:} Verification and validation procedures are used to compare the conceptual model with the proposed DES models. The development of the DES models is guided by the process flowchart and uses input data provided by GMH staff, who examined and approved these models. Additionally, the approach proposed by \Mycite{sargent2010verification} is adopted to verify and validate the DES models. \textbf{Data Validity:} The input data analysis section describes our data collection and analysis. This analysis indicates that our data is correct and adequately used. \textbf{Conceptual Model Validation:} The proposed conceptual model is validated via \emph{face validation} by GMH staff and via \emph{traces} following specific entities through the model. Flowcharts of the conceptual model are verified by GMH staff. Computerized Model Verification: The DES models are verified via techniques listed in \Mycite{sargent2010verification}. These techniques include animation, comparison with other models, and running several replications of the model. \textbf{Operational Validation:} A thorough sensitivity analysis was conducted to check the accuracy of the DES models. In the sensitivity analysis, the number of resources used (i.e., the number of AGVs, the number of instruments, etc.) changed, so the impact of these changes on the behavior of the model outputs was monitored. For example, the model outputs, after changing the number of AGVs, equal the average and standard deviation of travel time of the AGVs. Next, hypothesis tests were conducted to evaluate whether the difference between the outputs of DES models and the real-world data are statistically different. At a $p$-value of 0.05, the test indicates that the difference is not statistically significant.

\section{Discussion of Results}

The results from the DES models are used to address the research questions outlined in Section \ref{intro}.

\textbf{R1: How does the inventory level of surgical instruments, including owned, borrowed and cosigned, impact the efficiency of ORs?} A simulation-optimization experiment is conducted using ARENA Opt-Quest to answer this question. The objective of the simulation-optimization is to minimize the total delays at the start of a surgery by changing the inventory level. The delay of a surgery is calculated as the difference between the \textit{Actual Start Time} and the \textit{Scheduled Start Time}. The decision variables of type \textit{integer} are the number of instruments in the inventory for each of the seven service types (see Table \ref{Table: IISI}). In order to reduce the computational time of the simulation-optimization experiments, a lower bound, based on the data collected at GMH, is developed on the number of instruments inventoried. Let $n$ be the maximum number of surgeries scheduled in a day for each service type. The lower bound equals $\ceil{n/2}$. A lower bound is added for each surgery type via these constraints: (\textit{i}) the number of instruments used $\leq$ number of instruments in the inventory, and (\textit{ii}) the number of instruments used $\geq$ lower bound.

Experiments are conducted for three different scenarios. Scenario 1 assumes that the available inventory of instruments equals the current inventory level of GMH. Consider this inventory level to be an upper bound. Scenario 2 assumes that the available inventory of instruments equals the lower bound. Scenario 3 assumes that the available inventory of instruments equals the average value of the upper and lower bounds. Table \ref{Table:R1} summarizes the results of these experiments, and the following observations result:

\emph{Observation 1:} The \emph{Current} material handling approach, Model 1, is the most sensitive to changes in the inventory level, compared to \emph{Two Batch}, Model 2, and \emph{JIT}, Model 3. A decrease of inventory level, from Scenario 1 to Scenario 2, leads to an increase of the average delay from 0.42 to 31 minutes per surgery in the \emph{Current} approach. In the \emph{Two Batch} approach, the corresponding average delay increases from 0.01 to 5.12 minutes, and in the \emph{JIT} approach from 0.00 to 1.47 minutes per surgery (see Table \ref{Table:R1}).

\emph{Observation 2:} The \emph{Current} material handling approach requires additional levels of inventory to maintain the same service level, as measured by the average delay per surgery, compared to the proposed \emph{Two Batch} and \emph{JIT approaches} (see Table \ref{Table:R2}).

\emph{Observation 3:} The \emph{JIT} approach leads to reduced inventory levels of instruments used in short-duration surgeries without reducing the service level.

Table \ref{Table: NOS} presents the total number of neurological surgeries conducted at GMH during the 9-month period reviewed here. This number averages about 4.2 surgeries per day. Table \ref{Table: NOS} also presents the number of pediatric surgeries during the same time period, which corresponds to about 4.4 surgeries per day. The duration of neurological surgeries is about 1 hour longer than for pediatric surgeries. An hypothesis testing ($p$-value = 0.05) was conducted to evaluate the difference between the duration of neurological and pediatric surgeries. This test indicates that the difference is statistically significant (see Table \ref{tab:LengthComparison}). The results of Table \ref{Table:R2} show that the number of instruments required by neurological surgeries is higher than pediatric surgeries in Models 2 and 3 versus Model 1. This is because instruments used in Pediatric surgeries can be reused in the same day due to the shorter duration of these surgeries.

\begin{table}[H]
\caption{The Average Delay per Surgery}
\label{Table:R1}
\resizebox{\textwidth}{!}{%
\begin{tabular}{ccccccccccc}
\hline
\multicolumn{7}{c}{\textbf{Number of Instruments per Service Type}} & \textbf{} & \multicolumn{3}{c}{\textbf{Average Delay/Surgery (Minutes)}} \\ \hline
\multicolumn{1}{c|}{\textbf{Scenario}} & \textbf{ENT} & \textbf{Gynecology} & \textbf{Neurological} & \textbf{Ortho Trauma} & \textbf{Pediatric} & \textbf{Urology} & \multicolumn{1}{c|}{\textbf{Vascular}} & \textbf{Model 1} & \textbf{Model 2} & \textbf{Model 3} \\ \hline
\multicolumn{1}{c|}{\textbf{1}} & 10 & 10 & 13 & 17 & 8 & 18 & \multicolumn{1}{c|}{16} & 0.42 & 0.01 & 0 \\
\multicolumn{1}{c|}{\textbf{2}} & 6 & 6 & 5 & 9 & 4 & 10 & \multicolumn{1}{c|}{6} & 31.27 & 5.12 & 1.47 \\
\multicolumn{1}{c|}{\textbf{3}} & 8 & 8 & 9 & 13 & 6 & 14 & \multicolumn{1}{c|}{10} & 3.58 & 0.25 & 0.01 \\
\hline
\end{tabular}%
}
\end{table}

\begin{table}[H]
\caption{Inventory Level of Instruments}
\label{Table:R2}
\resizebox{\textwidth}{!}{%
\begin{tabular}{ccccccccc}
\hline
 &  & \multicolumn{7}{c}{\textbf{Number of Instruments per Service}} \\ \hline
\multicolumn{1}{c|}{\textbf{Model}} & \multicolumn{1}{c|}{\textbf{Delay/ Surgery (Minutes)}} & \textbf{ENT} & \textbf{Gynecology} & \textbf{Neurological} & \textbf{Ortho Trauma} & \textbf{Pediatric} & \textbf{Urology} & \textbf{Vascular} \\ \hline
\multicolumn{1}{c|}{\textbf{1}} & \multicolumn{1}{c|}{0.42} & 10 & 10 & 13 & 17 & 8 & 18 & 16 \\
\multicolumn{1}{c|}{\textbf{2}} & \multicolumn{1}{c|}{0.41} & 6 & 10 & 13 & 13 & 6 & 18 & 12 \\
\multicolumn{1}{c|}{\textbf{3}} & \multicolumn{1}{c|}{0.41} & 6 & 6 & 9 & 9 & 4 & 16 & 12 \\ \hline
\end{tabular}%
}
\end{table}

\begin{table}[H]
\caption{Comparison of Two Service Types}
\label{tab:LengthComparison}
\centering
\resizebox{0.55\textwidth}{!}{%
\begin{tabular}{l|cc}
\hline
\textbf{Statistics} & \textbf{Neurology Surgery} & \textbf{Pediatric Surgery} \\ \hline
\textbf{Sample Size} & 1062 & 1121 \\
\textbf{Average Length} & 2.31 & 1.36 \\
\textbf{95\% CI} & (2.24,2.38) & (1.31,1.40) \\
\textbf{Standard Deviation} & 1.16 & 0.77 \\ \hline
\end{tabular}%
}

\end{table}

\textbf{R2: How do material handling activities impact the efficiency of ORs?} Two sets of experiments are conducted. The first set focuses on the impact that changing the number of AGVs has on the performance of the material handling system. This performance is measured via the average travel time per trip, the total travel time, and the corresponding standard deviations. The delivery time of clean and soiled case carts is analyzed as the number of AGVs increases from 6 to 8 to 10. Experiments with fewer than 6 AGVs led to extensive delays in delivering all the case carts in \emph{Current} system, which requires employees to work overtime, so these experiments are not considered in this analysis.

\begin{figure}[H]
    \centering
    \begin{subfigure}[b]{0.3\textwidth}
        \includegraphics[width=\textwidth]{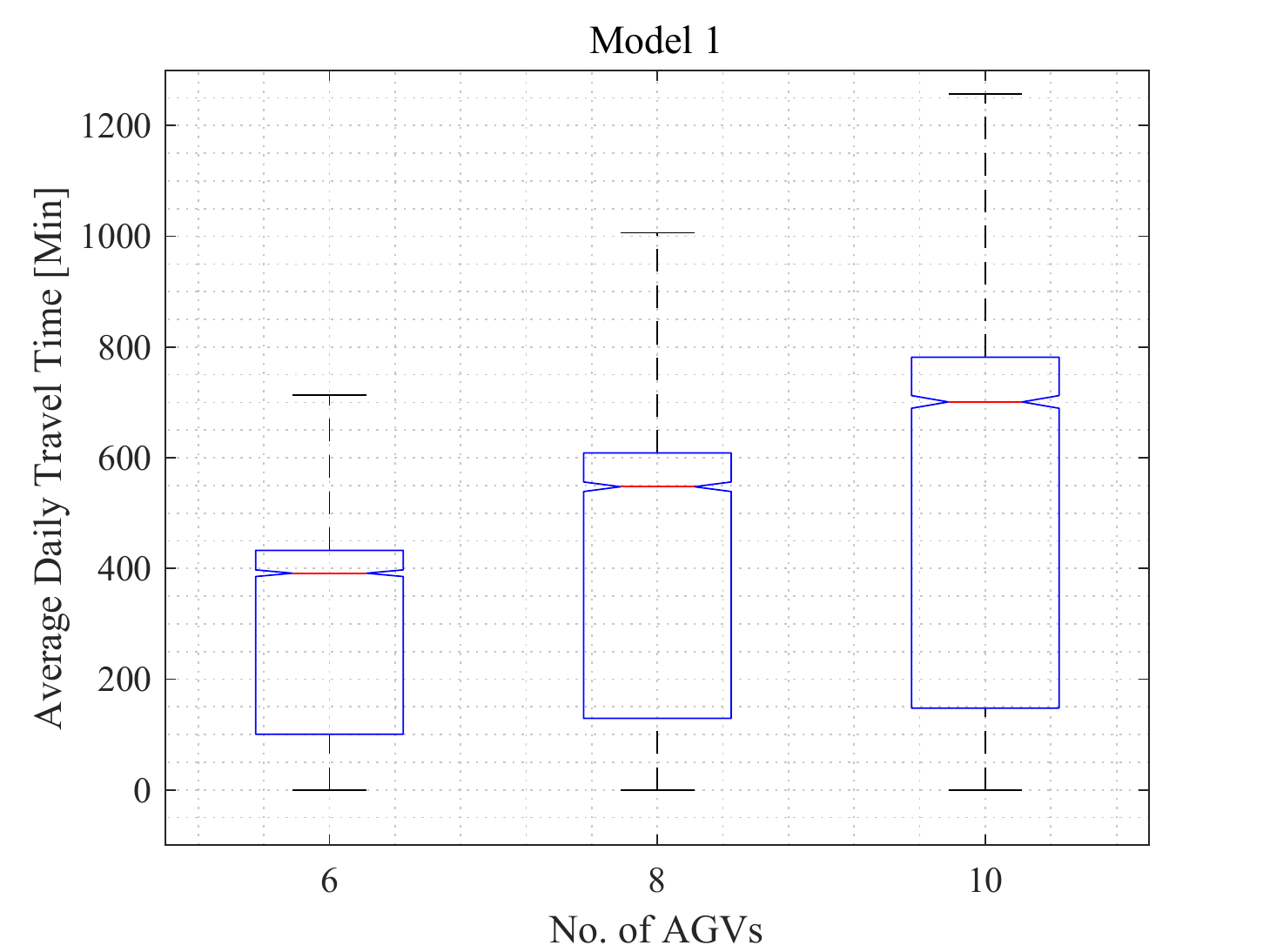}
        \caption{Total Travel Time: Model 1}
        \label{fig:cm1}
    \end{subfigure}
    ~ 
    \begin{subfigure}[b]{0.3\textwidth}
        \includegraphics[width=\textwidth]{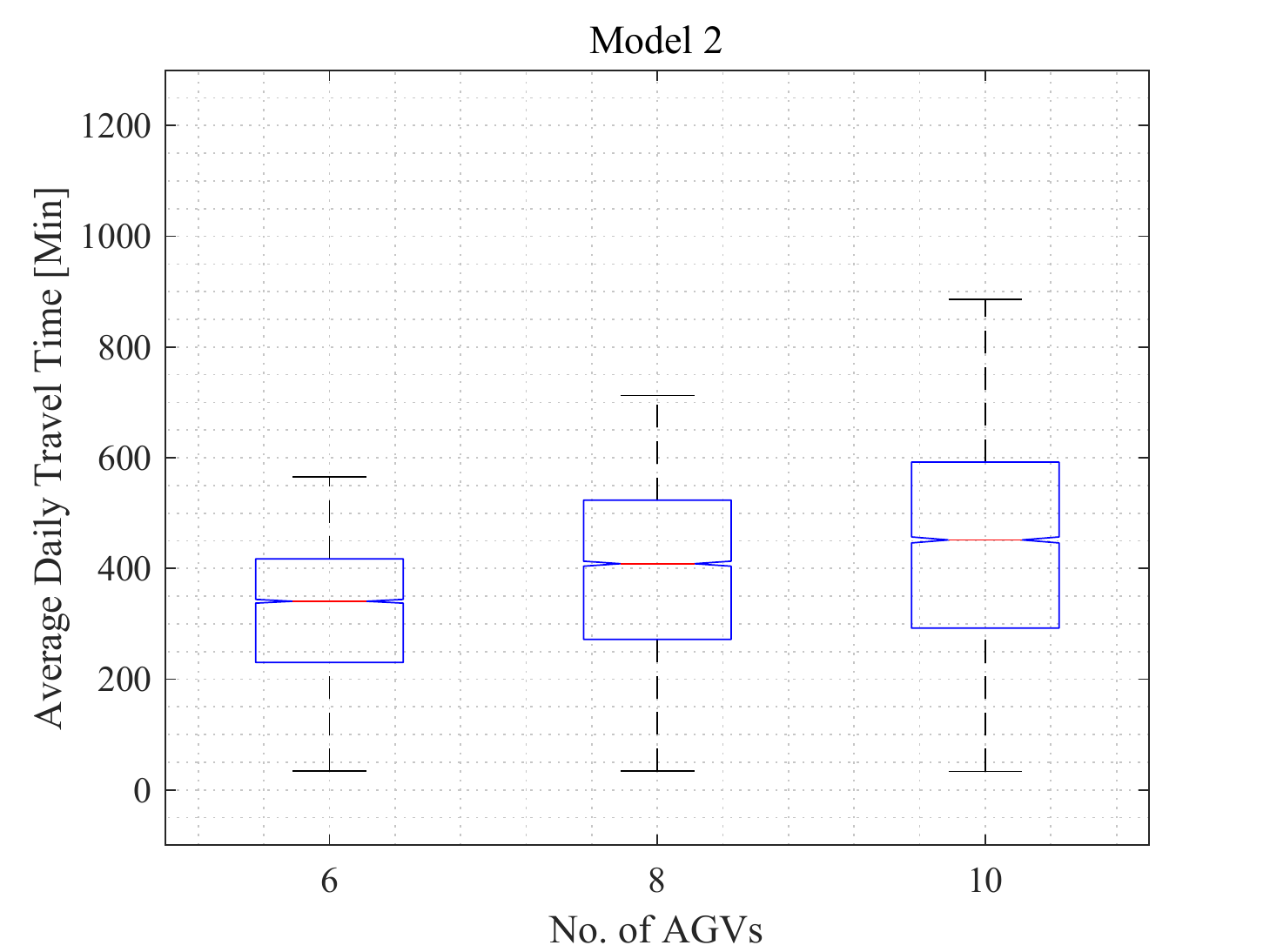}
        \caption{Total Travel Time: Model 2}
        \label{fig:cm2}
    \end{subfigure}
    ~ 
     \begin{subfigure}[b]{0.3\textwidth}
        \includegraphics[width=\textwidth]{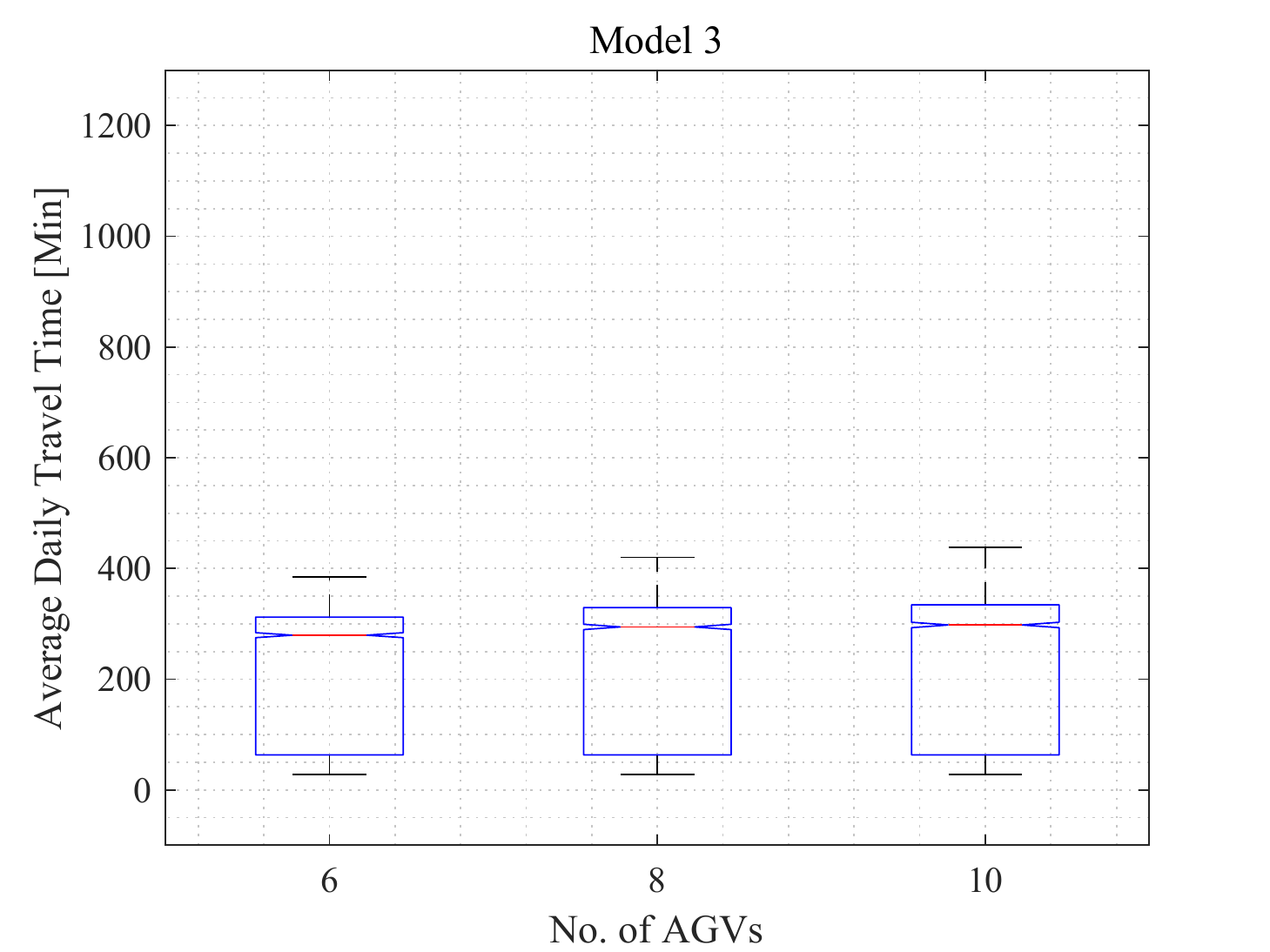}
        \caption{Total Travel Time: Model 3}
        \label{fig:cm3}
    \end{subfigure}
    \caption{Sensitivity Analysis of Clean Case Carts Delivery}
    \label{fig:EffectofMaterialHandlingCCS}
\end{figure}

\begin{figure}[ht]
    \centering
    \begin{subfigure}[b]{0.3\textwidth}
        \includegraphics[width=\textwidth]{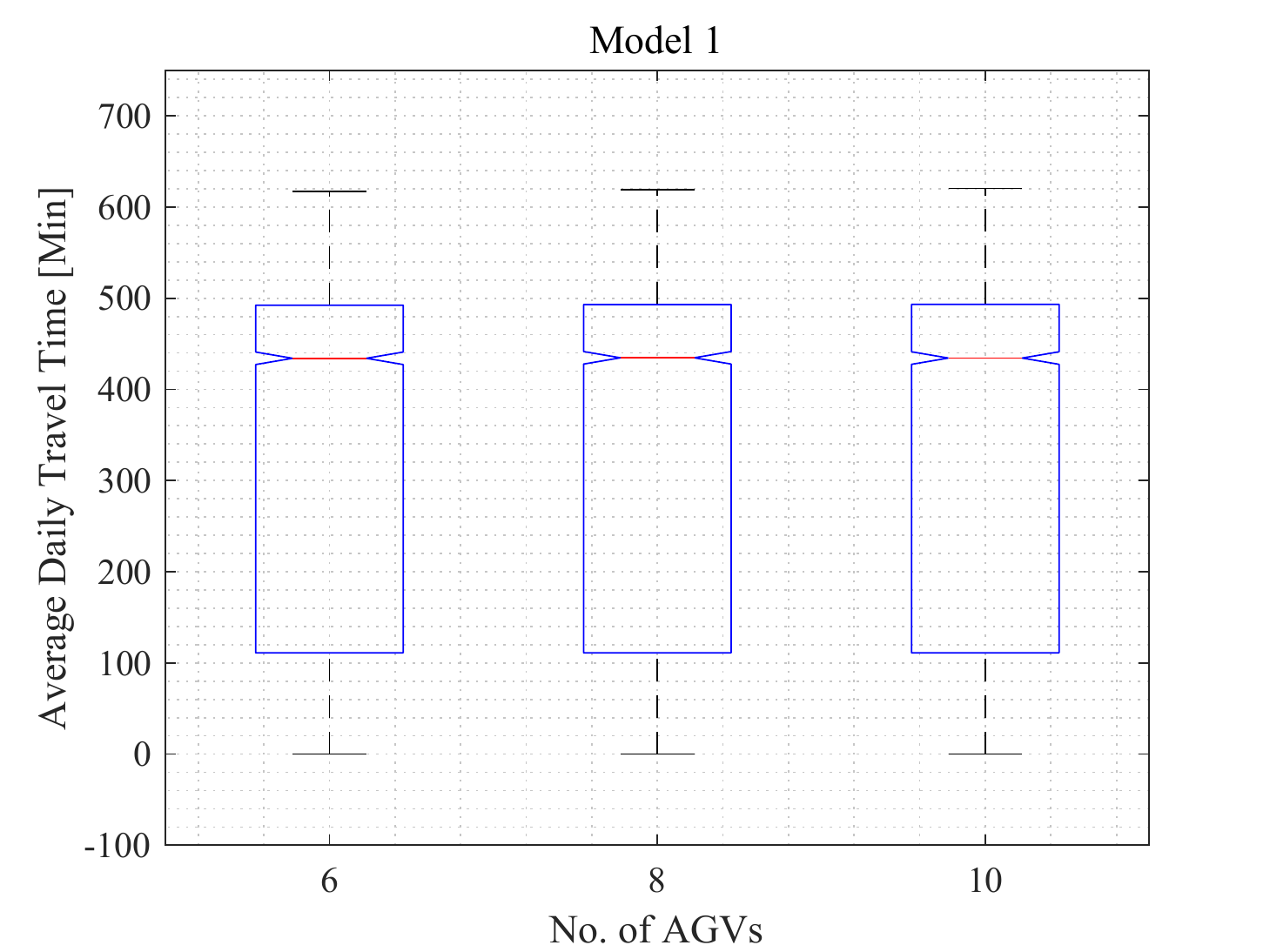}
        \caption{Total Travel Time: Model 1}
        \label{fig:dm1}
    \end{subfigure}
    ~ 
    \begin{subfigure}[b]{0.3\textwidth}
        \includegraphics[width=\textwidth]{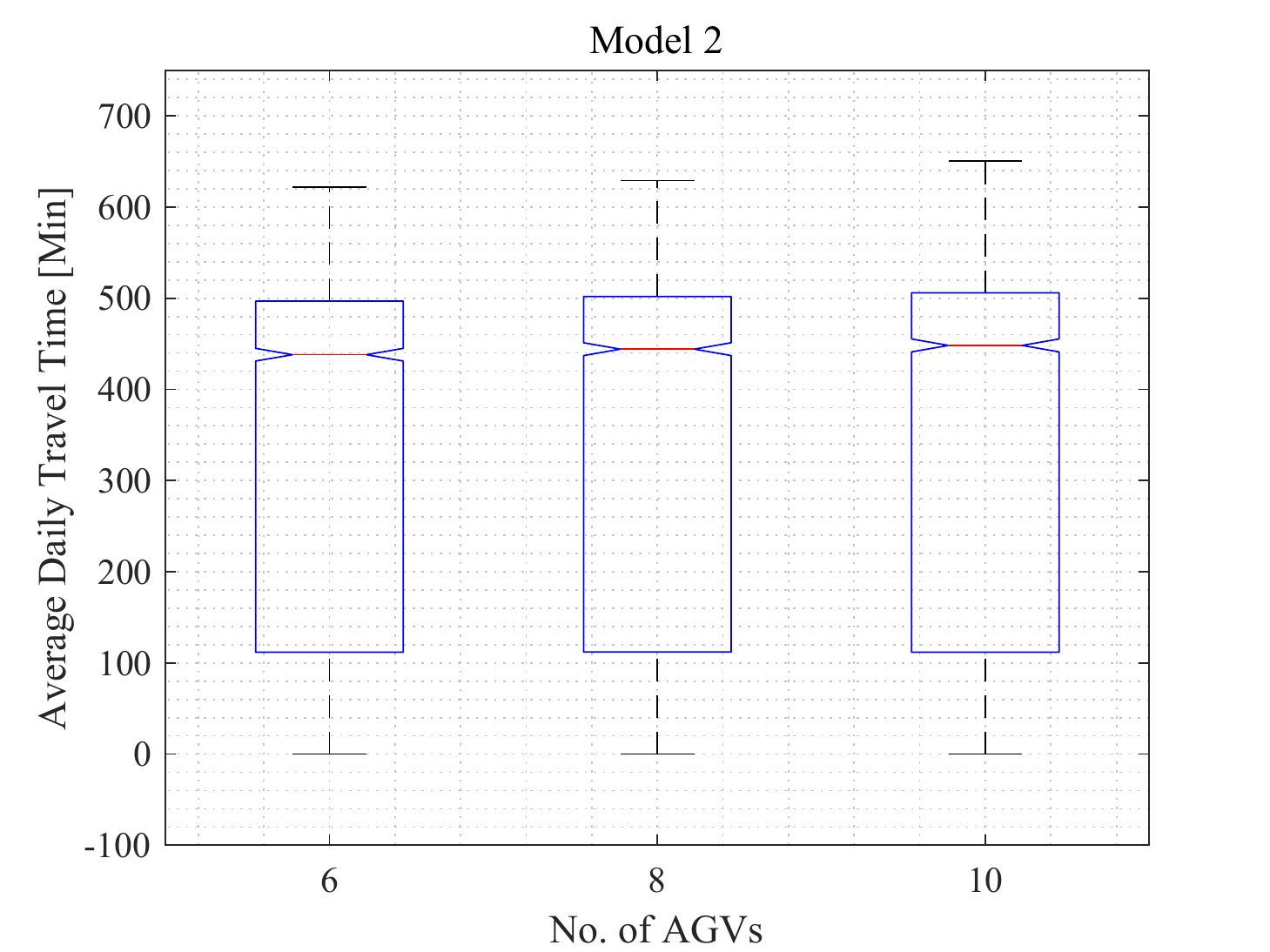}
        \caption{Total Travel Time: Model 2}
        \label{fig:dm2}
    \end{subfigure}
    ~ 
     \begin{subfigure}[b]{0.3\textwidth}
        \includegraphics[width=\textwidth]{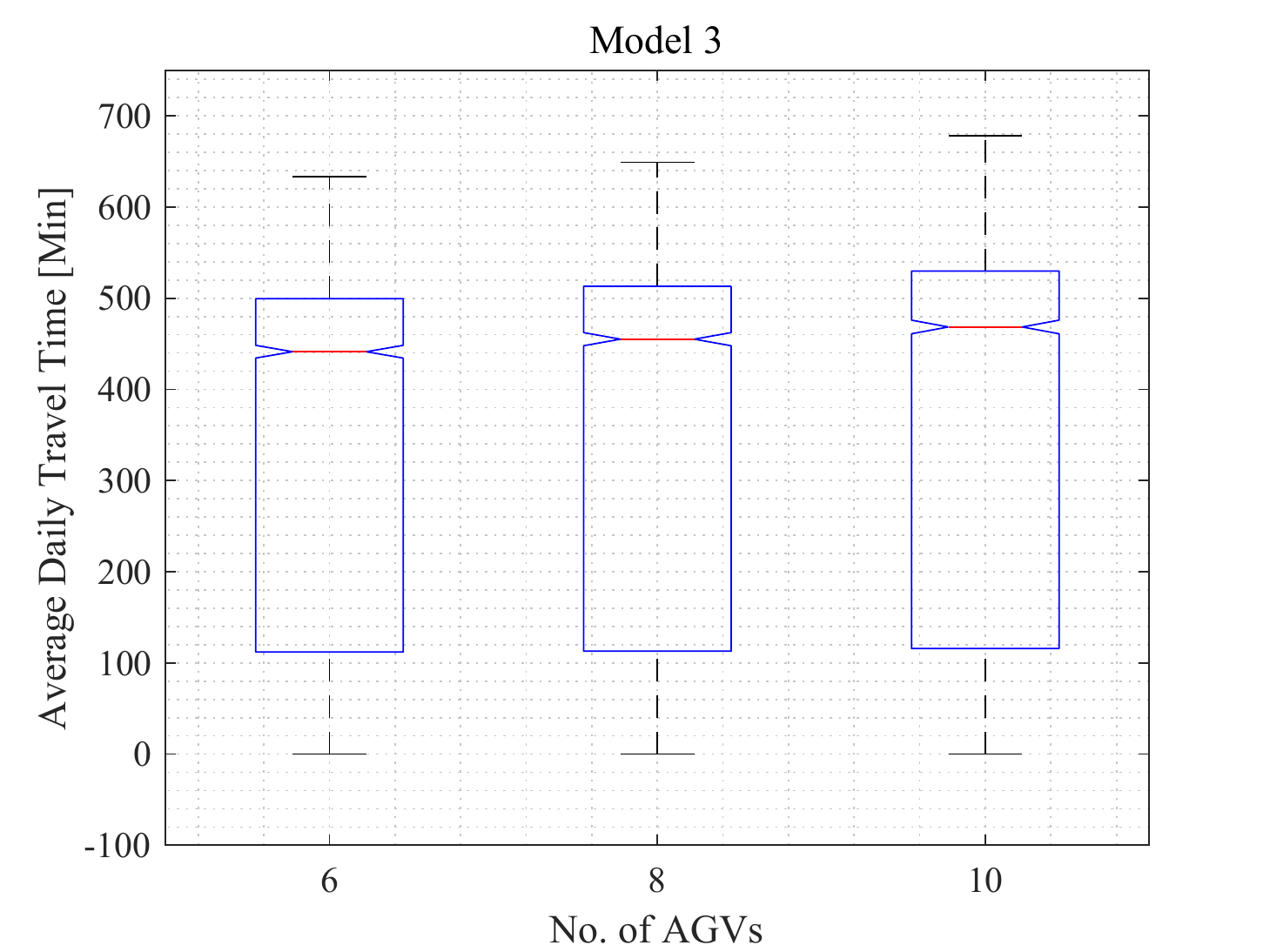}
        \caption{Total Travel Time: Model 3}
        \label{fig:dm3}
    \end{subfigure}
    \caption{Sensitivity Analysis of Soiled Case Carts Delivery}
    \label{fig:EffectofMaterialHandlingSCS}
\end{figure}

The second set of experiments focuses on the impact that changing the timing of delivery has on the performance of the material handling system. For this purpose, the performances of Models 1, 2, and 3 are compared. The results of these experiments are summarized in Tables \ref{Table:CCS} and \ref{Table:SCS} and Figures \ref{fig:EffectofMaterialHandlingCCS} and \ref{fig:EffectofMaterialHandlingSCS}. The following observations result:

\emph{Observation 1:} The average daily travel time of clean case carts is longest in the \emph{Current} material handling approach and shortest in the \emph{JIT} approach (see Figures \ref{fig:cm1}, \ref{fig:cm2}, and \ref{fig:cm3}). The \emph{Current} approach has the longest travel time due to congestion since the delivery of clean case carts for elective surgeries takes place during 3-7pm.

\emph{Observation 2:} The average daily travel time of clean case carts increases with the number of AGVs (see Figures \ref{fig:cm1}, \ref{fig:cm2}, and \ref{fig:cm3}). This increase is highest in the \emph{Current} material handling approach.

\emph{Observation 3:} The average daily travel time of clean case carts in the \emph{JIT} approach is not impacted by the increase in the number of AGVs since the delivery of case carts is spread over the day. These deliveries do not cause congestion (see Figure \ref{fig:cm3}).

\emph{Observation 4:} The average daily travel time of soiled case carts for every material handling approach is slightly impacted by the increase in the number of AGVs (see Figures \ref{fig:dm1}, \ref{fig:dm2}, and \ref{fig:dm3}). Note that the difference in the average travel time per trip is small but still statistically significant. The change in travel time due to the increase in the number of AGVs for every Model is small because soiled case carts are delivered to CSSD right after the surgery; thus, they are delivered throughout the day, and these deliveries have a minimal impact on congestion.

\begin{table}[ht]
\caption{Sensitivity Analysis of Clean Case Carts Delivery}
\label{Table:CCS}
\centering
\resizebox{0.7\textwidth}{!}{%
\begin{tabular}{c|c|cccc}
\hline
\textbf{} & \textbf{} & \multicolumn{4}{c}{\textbf{Travel Time (Minutes)}} \\ \hline
\textbf{No. of AGVs} & \textbf{Model} & \textbf{Average} & \textbf{StDev} & \textbf{CI for Average} & \textbf{CI for StDev} \\ \hline
\textbf{6} & \textbf{1} & 4.83 & 0.31 & (4.83, 4.84) & (0.30, 0.31) \\
\textbf{} & \textbf{2} & 4.96 & 0.17 & (4.95, 4.96) & (0.17, 0.18) \\
\textbf{} & \textbf{3} & 3.31 & 0.27 & (3.31, 3.32) & (0.26, 0.27) \\ \hline
\textbf{8} & \textbf{1} & 6.53 & 0.60 & (6.51, 6.54) & (0.59, 0.61) \\
\textbf{} & \textbf{2} & 6.47 & 0.64 & (6.46, 6.48) & (0.62, 0.65) \\
\textbf{} & \textbf{3} & 3.43 & 0.36 & (3.42, 3.44) & (0.35, 0.36) \\ \hline
\textbf{10} & \textbf{1} & 8.02 & 1.16 & (7.99, 8.05) & (1.14, 1.17) \\
\textbf{} & \textbf{2} & 7.66 & 1.17 & (7.63, 7.69) & (1.14, 1.19) \\
\textbf{} & \textbf{3} & 3.46 & 0.39 & (3.46, 3.47) & (0.38, 0.39)
\\ \hline
\end{tabular}%
}

\end{table} 

\begin{table}[ht]
\centering
\caption{Sensitivity Analysis of Soiled Case Carts Delivery}
\label{Table:SCS}
\resizebox{0.7\textwidth}{!}{%
\begin{tabular}{c|c|cccc}
\hline
\textbf{} & \textbf{} & \multicolumn{4}{c}{\textbf{Travel Time (Minutes)}} \\ \hline
\textbf{No. of AGVs} & \textbf{Model} & \textbf{Average} & \textbf{StDev} & \textbf{CI for Average} & \textbf{CI for StDev} \\ \hline
\textbf{6} & \textbf{1} & 5.46 & 0.10 & (5.46, 5.46) & (0.094, 0.097) \\
\textbf{} & \textbf{2} & 5.44 & 0.09 & (5.44, 5.44) & (0.087, 0.089) \\
\textbf{} & \textbf{3} & 5.39 & 0.07 & (5.39, 5.39) & (0.069, 0.071) \\ \hline
\textbf{8} & \textbf{1} & 5.60 & 0.18 & (5.59, 5.60) & (0.179, 0.186) \\
\textbf{} & \textbf{2} & 5.53 & 0.15 & (5.53, 5.54) & (0.149, 0.154) \\
\textbf{} & \textbf{3} & 5.40 & 0.07 & (5.39, 5.40) & (0.072, 0.074) \\ \hline
\textbf{10} & \textbf{1} & 5.75 & 0.29 & (5.75, 5.76) & (0.290, 0.299) \\
\textbf{} & \textbf{2} & 5.64 & 0.23 & (5.63, 5.64) & (0.232, 0.238) \\
\textbf{} & \textbf{3} & 5.40 & 0.07 & (5.39, 5.40) & (0.072, 0.074) \\ \hline
\end{tabular}%
}
\end{table}

\textbf{R3: How does integrating decisions about inventory and material handling impact the efficiency of ORs?}

\emph{Observation 1:} The average delay per surgery and the total number of delays are lowest in the \emph{JIT} material handling approach. These statistics are highest in the \emph{Current} approach. A successful implementation of \emph{JIT} requires coordination of material handling and inventory decisions, and numerical results show that this coordination leads to improved OR efficiency.

This observation is true at different inventory levels, represented by Scenarios 1, 2, and 3 in Table \ref{Table:R1}; for different material handling approaches, represented by Models 1, 2, and 3 in Tables \ref{tab:AvgDelayByType} and \ref{tab:delayfrequency}; and for different material handling capacities, represented by the number of AGVs in Tables \ref{tab:AvgDelayByType} and \ref{tab:delayfrequency}.

\emph{Observation 2:} The average inventory level is lowest in the \emph{JIT} material handling approach. This observation is supported by the results of Tables \ref{Table:R1} and \ref{Table:R2}.

\textbf{Recommendations:} The following recommendations are made based on the observations presented above:

\emph{Recommendation 1:} Coordinating material handling and inventory management decisions has the potential to improve the efficiency and reduce the cost of ORs without jeopardizing the service level provided. To facilitate this coordination, the requirements for each surgery, the number of available instruments, and the location of instruments must be known at all times. Transparent information technology systems will facilitate the coordination of decisions.

\emph{Recommendation 2:}  Hospitals should consider implementing a \emph{JIT} material handling approach for instruments used in short-duration surgeries because such an approach leads to lower inventory levels without jeopardizing the service level provided.  

\emph{Recommendation 3:} Hospitals should frequently reevaluate their material handling system to identify improvements. For example, GMH currently uses 10 AGVs. Using only 6 or 8 AGVs leads to reduced congestion along the corridors of the hospital and leads to shorter delivery times. The remaining AGVs can be used for transportation of trash, linen, and pharmaceuticals, among other items.

\begin{table}[H]
\caption{Average Delay by Service Type (Hours)}
\label{tab:AvgDelayByType}
\resizebox{\textwidth}{!}{%
\begin{tabular}{c|ccc|ccc|ccc}
\hline
\textbf{} & \multicolumn{3}{c|}{\textbf{Model 1}} & \multicolumn{3}{c|}{\textbf{Model 2}} & \multicolumn{3}{c}{\textbf{Model 3}} \\ \hline
\textbf{Service Type} & \textbf{6 AGV} & \textbf{8 AGV} & \textbf{10 AGV} & \textbf{6 AGV} & \textbf{8 AGV} & \textbf{10 AGV} & \textbf{6 AGV} & \textbf{8 AGV} & \textbf{10 AGV} \\ \hline
ENT  & 0.060 & 0.006 & 0.001 & 0.009 & 0.000 & 0.000 & 0.001 & 0.000 & 0.000 \\
Gynecology  & 0.132 & 0.011 & 0.001 & 0.014 & 0.000 & 0.000 & 0.000 & 0.000 & 0.000 \\
Neurological  & 1.672 & 0.124 & 0.005 & 0.330 & 0.011 & 0.000 & 0.085 & 0.000 & 0.000 \\
Ortho trauma  & 0.018 & 0.000 & 0.000 & 0.001 & 0.000 & 0.000 & 0.000 & 0.000 & 0.000 \\
Pediatric  & 1.186 & 0.244 & 0.042 & 0.266 & 0.020 & 0.002 & 0.030 & 0.001 & 0.000 \\
Urology  & 0.245 & 0.024 & 0.001 & 0.014 & 0.000 & 0.000 & 0.000 & 0.000 & 0.000 \\
Vascular  & 0.324 & 0.003 & 0.000 & 0.028 & 0.000 & 0.000 & 0.003 & 0.000 & 0.000 \\ \hline
\multicolumn{3}{c}{\footnotesize{*The total number of replications is 30.}}
\end{tabular}%
}
\end{table}

\begin{table}[H]
\caption{Frequency of Delayed Surgeries}
\resizebox{\textwidth}{!}{%
\begin{tabular}{c|ccc|ccc|ccc}
\hline
 & \multicolumn{3}{c|}{\textbf{Model 1}} & \multicolumn{3}{c|}{\textbf{Model 2}} & \multicolumn{3}{c|}{\textbf{Model 3}} \\ \hline
\textbf{Service Type} & \textbf{6 AGV} & \textbf{8 AGV} & \textbf{10 AGV} & \textbf{6 AGV} & \textbf{8 AGV} & \textbf{10 AGV} & \textbf{6 AGV} & \textbf{8 AGV} & \textbf{10 AGV} \\ \hline
{ENT} & 378 & 41 & 5 & 175 & 14 & 2 & 17 & 0 & 0 \\
{Gynecology} & 580 & 60 & 5 & 138 & 10 & 0 & 10 & 0 & 0 \\
{Neurological} & 6,254 & 595 & 22 & 3,352 & 314 & 0 & 1,820 & 0 & 0 \\
{Ortho Trauma} & 120 & 0 & 0 & 63 & 0 & 0 & 0 & 0 & 0 \\
{Pediatric} & 4,859 & 1,293 & 255 & 2,532 & 377 & 52 & 743 & 20 & 0 \\
{Urology} & 2,261 & 268 & 8 & 652 & 39 & 0 & 8 & 0 & 0 \\
{Vascular} & 1,722 & 24 & 0 & 615 & 5 & 0 & 85 & 0 & 0 \\ \hline
\multicolumn{3}{c}{\footnotesize{*The total number of replications is 30.}}
\end{tabular}%
}
\label{tab:delayfrequency}
\end{table}

\section{Summary and Concluding Remarks}

The proposed research and the models presented here are motivated by the opportunities for improvement observed in GMH’s inventory management and material handling. Based on their \emph{Current} material handling system, case carts loaded with instruments required by surgeries scheduled during a day are delivered, via AGVs, to ORs within a few hours or the evening before the surgery. This delivery schedule leads to (\textit{i}) increased inventory of instruments, whether owned, loaned, or cosigned; (\textit{ii}) increased traffic and congestion, which delay the deliveries of case carts and the delivery other materials that use AGVs; and (\textit{iii}) delayed surgery start times. The inefficiencies identified motivated the following research questions: \textbf{(R1)}: \emph{How does the inventory level of surgical instruments, including owned, borrowed and cosigned, impact the efficiency of ORs?} \textbf{(R2)}: \emph{How do material handling activities impact the efficiency of ORs?} \textbf{(R3)}: \emph{How does integrating decisions about inventory and material handling impact the efficiency of ORs?}

In order to address these research questions, two new material handling approaches are proposed and compared to the approach already in place. Along with the already-existing \emph{Current} approach, the \emph{Two Batch} approach delivers surgical carts to ORs twice a day, in the morning and in the evening, and the \emph{JIT} approach delivers a surgical carts to an OR before the surgery. Three DES models are developed for each of the three approaches, and they are verified and validated using real-life data collected at a partnering hospital.

A thorough sensitivity analysis of the DES models is conducted and leads to a number of observations and recommendations. The \emph{Current} material handling approach is most sensitive to changes in the inventory level, requires the highest levels of inventory to maintain a high service level, and leads to congestion and delays of the delivery of surgical case carts. Both the \emph{Two Batch} and \emph{JIT} approaches outperform the Current material handling approach. The implementation of the \emph{JIT} approach leads to the greatest improvements in OR efficiency and service levels. Based on these observations, hospitals should identify opportunities to coordinate material handling and inventory management decisions since it leads to improved efficiency for ORs. New, data-based approaches to material handling and inventory management, like the \emph{JIT} delivery of surgical cases, have the potential to improve the efficiency of short-duration surgeries in hospital ORs.

\subsection*{ACKNOWLEDGMENTS} This research is partially funded via a Spark Research Grant awarded by Material Handling Institute (MHI) and College-Industry Council on Material Handling Education (CICMHE). The authors are grateful to HMI and CICMHE for their support. The authors are thankful to the staff of Greenville Memorial Hospital for their continuous support in this research by providing the data, and the expertise necessary to verify and validate the models developed.

\bibliographystyle{plainnat}
\bibliography{FinalPaper}
\end{document}